\documentstyle[11pt,twoside]{article} 
\def\be{\begin{equation}}
\def\ee{\end{equation}}
\def\bea{\begin{eqnarray}}
\def\eea{\end{eqnarray}}
\def\nn{\nonumber}
\def\ba{\begin{array}}
\def\ea{\end{array}}
\def\lb{\left(}
\def\rb{\right)}
\def\lsb{\left[}
\def\rsb{\right]}
\def\lc{\left\{}

\def\s{\sum_{n=0}^\infty}
\def\a{\alpha}
\def\b{\beta}
\def\t{\tilde}

\begin{document}
\begin{flushright}
math/9803142
\end{flushright}
\begin{flushright}
IMSc-98/03/11
\footnote{{\sf Published in} {\em Special Functions and Differential 
Equations}, Proceedings of a Workshop held at The Institute of 
Mathematical Sciences, Madras, India, January 13-24, 1997, Eds. 
K. Srinivasa Rao, R. Jagannathan, G. Vanden Berghe and J. Van der Jeugt 
(Allied Publishers, New Delhi, 1998) pp. 158-164.}
\end{flushright} 
\begin{center} {\Large \bf $(P,Q)$-Special Functions} 

\medskip

{\bf R. Jagannathan} \\ 
The Institute of Mathematical Sciences \\ 
C.I.T. Campus, Tharamani, Chennai~(Madras), Tamilnadu - 600~113, India. \\
e-mail:~{\tt jagan@imsc.ernet.in}
\end{center}

\bigskip

\noindent{\bf Abstract:}\ \ It is suggested that the 
$(p,q)$-hypergeometric series studied by Burban and Klimyk (in 
{\em Integral Transforms and Special Functions} {\bf 2} (1994) 
15 - 36) can be considered as a special case of a more general 
$(P,Q)$-hypergeometric series. 

\vspace{.5cm}

\noindent{\bf 1. Introduction}

\smallskip

\noindent I propose a general $(P,Q)$-hypergeometric series.  
In this, I am inspired mainly by the paper of Burban and Klimyk~\cite{BK} 
titled ``$P,Q$-differentiation, $P,Q$-integration, and 
$P,Q$-hypergeometric functions related to quantum groups'' and 
the paper of Floreanini, Lapointe and Vinet~\cite{FLV} titled ``A 
note on $(p,q)$-oscillators and bibasic hypergeometric functions''.  
Burban and Klimyk~\cite{BK} have already presented a well-developed 
theory of $(p,q)$-hypergeometric functions and what I propose is only a 
suggestion towards a slight generalization of their work, taking some 
clues from Floreanini, Lapointe and Vinet~\cite{FLV}, Katriel and 
Kibler~\cite{KK}, and Gasper and Rahman~\cite{GR}.  

\bigskip

\noindent{\bf 2. $q$-Hypergeometric Series} 
\renewcommand{\theequation}{2.{\arabic{equation}}}
\setcounter{equation}{0}

\smallskip

\noindent One defines the $q$-shifted factorial by 
\be
(x;q)_n = \lc 
\ba{ll}
1, & n = 0, \\
(1-x)(1-xq)(1-xq^2) \dots \\ 
\quad \dots\,\lb 1-xq^{n-1} \rb , & n = 1,2,\dots . 
\ea \right. 
\ee
Then, with the notation
\be
\lb x_1,x_2,\,\dots\,,x_m ; q \rb_n = 
\lb x_1;q \rb _n \lb x_2;q \rb _n\,\dots\,\lb x_m;q \rb _n \,, 
\ee
an $_r\phi_s$ basic hypergeometric series, or a general 
$q$-hypergeometric series, is given by  
\bea
   &   & 
_r\phi_s \lb a_1, a_2,\,\dots\,,a_r;b_1,b_2,\,\dots\,,b_s;q,z \rb \nn \\ 
   &   & \quad = \s \,\frac{ \lb a_1,a_2,\,\dots\,,a_r ; q \rb_n }
{ \lb b_1,b_2,\,\dots\,,b_s ; q \rb_n (q;q)_n } 
\lb (-1)^n q^{n(n-1)/2} \rb ^{1+s-r} z^n\,,
\label{phirs}
\eea
with $|q| < 1$ and $r,s = 0,1,2,\,\ldots $ (see~\cite{GR} for details).  
Choosing 
\be
a_1 = q^{\a _1},\ \ \dots\,,\ \ a_r = q^{\a _r}, \quad    
b_1 = q^{\b _1},\ \ \dots\,,\ \ b_s = q^{\b _s}, 
\ee 
and defining, following Heine, 
\bea
[x]_q & = & \frac{1-q^x}{1-q}\,, \\
   &    &  \nn \\ 
\lb [x]_q \rb _n & = & \lc 
\ba{ll}
1, & n = 0, \cr  
[x]_q [x+1]_q [x+2]_q \dots & \cr 
\qquad \dots\,[x+n-1]_q, & n = 1,2,\dots , 
\ea \right. \\ 
   &    &  \nn \\ 
{[n]_q}! & = & \lc 
\ba{ll}
1, & n = 0, \cr
[n]_q [n-1]_q [n-2]_q \dots &  \cr  
\qquad \dots\,[2]_q [1]_q, & n = 1,2,\dots , 
\ea \right.
\eea
one gets the special case: 
\bea
   &   & 
_r\phi_s \lb q^{\a _1},q^{\a_ 2},\dots ,q^{\a_ r};q^{\b _1},q^{\b _2},
\dots ,q^{\b _s};q,z \rb \nn \\ 
   &   & \quad = \s \,\frac{ \lb \lsb \a _1 \rsb _q \rb _n 
\lb \lsb \a _2 \rsb _q \rb _n \dots \lb \lsb \a _r \rsb _q \rb _n }
{ \lb \lsb \b _1 \rsb _q \rb _n \lb \lsb \b _2 \rsb _q \rb _n 
\dots \lb \lsb \b _s \rsb _q \rb _n }\, 
\frac{ \lb (-1)^n q^{n(n-1)/2} \rb ^{(1+s-r)} }{ (1-q)^{n(1+s-r)} 
[n]_q! } \,z^n \,. 
\eea  

\bigskip

\noindent{\bf 3. $(P,Q)$-Hypergeometric Series}
\renewcommand{\theequation}{3.{\arabic{equation}}}
\setcounter{equation}{0}

\smallskip

\noindent Now comes my proposal.  Let  
\be
\t{q} = (P,Q) 
\ee 
and 
\bea
((P,Q);(P,Q))_n & = & \lb \t{q};\t{q} \rb _n \nn \\ 
   & = & \lc 
\ba{ll}
1, & n = 0 \cr 
(P-Q) \lb P^2-Q^2 \rb \lb P^3-Q^3 \rb \,\dots & \cr 
\qquad \dots\,\lb P^n-Q^n \rb \,, & n = 1,2,\,\dots \,. 
\ea \right. 
\label{qqn}
\eea
For any 
\be
\t{x} = \lb x_p,x_q \rb \,, 
\ee 
let
\bea
\lb \lb x_p,x_q \rb ; (P,Q) \rb _n & = & 
\lb \t{x};\t{q} \rb _n \nn \\ 
   & = & \lc 
\ba{ll}
1, & n = 0, \cr
\lb x_p-x_q \rb \lb x_pP-x_qQ \rb \lb x_pP^2-x_qQ^2 \rb \dots & \cr 
\qquad \dots\,\lb x_pP^{n-1}-x_qQ^{n-1} \rb , &  
n = 1,2,\,\dots\,. 
\ea \right. 
\label{xqn}
\eea 
As before, I shall use the notation 
\be
\lb \t{x}_1,\t{x}_2,\,\dots\,,\t{x}_m ; \t{q} \rb_n =  
\lb \t{x}_1; \t{q} \rb _n \lb \t{x}_2; \t{q} \rb _n\,\dots\,\lb 
\t{x}_m; \t{q} \rb _n \,.
\ee
Then, with  
\bea
   &   & \t{a}_1 = \lb a_{1p},a_{1q} \rb , \ \ \dots \,,\ \ 
\t{a}_r = \lb a_{rp},a_{rq} \rb \,, \nn \\ 
   &   & \t{b}_1 = \lb b_{1p},b_{1q} \rb , \ \ \dots \,,\ \ 
\t{b}_s = \lb b_{sp},b_{sq} \rb \,,  
\eea
I define an $_r\tilde{\phi}_s$ basic hypergeometric series, or a 
general $(P,Q)$-hypergeometric series, by 
\bea
   &   & 
_r\tilde{\phi}_s \lb \t{a}_1,\t{a}_2,\dots ,
\t{a}_r;\t{b}_1,\t{b}_2, \dots ,\t{b}_s;\t{q},z \rb \nn \\ 
   &   & \ = \s \,\frac{ \lb \t{a}_1,\t{a}_2,\,\dots\,,\t{a}_r; 
\t{q} \rb_n } { \lb \t{b}_1,\t{b}_2,\,\dots\,,\t{b}_s; \t{q} \rb_n 
\lb \t{q};\t{q} \rb _n } \lb (-1)^n (Q/P)^{n(n-1)/2} \rb 
^{1+s-r} z^n\,. 
\label{tphi}
\eea
with $|Q/P| < 1$ and $r,s = 0,1,2,\,\ldots$\,\,. \\

The $_r\phi_s$ series~(\ref{phirs}) is a special case of 
$_r\t{\phi}_s$ series~(\ref{tphi}) corresponding to 
the choice $a_{1p} = a_{2p} =\,\dots\,= a_{rp} = b_{1p} = b_{2p} 
=\,\dots\,= b_{sp} = 1$, $a_{1q} = a_1$, $a_{2q} = a_2$,\  $\ldots$\,,
\  $a_{rq} = a_r$, $b_{1q} = b_1$, $b_{2q} = b_2$,\  $\ldots$\,,\   
$b_{sq} = b_s$ and $(P,Q) = (1,q)$.  The factor 
$\lb (-1)^n q^{n(n-1)/2} \rb ^{1+s-r}$ included in the 
definition of $_r\phi_s$ series~(\ref{phirs}), following Gasper and 
Rahman~\cite{GR}, is absent in the earlier literature (see, 
{\em e.g.},~\cite{B,A,S}).  The inclusion of this factor leads to 
the nice property 
\be
\lim_{a_r \rightarrow \infty}\,_r\phi_s \lb z/a_r \rb 
=\,_{r-1}\phi_s (z)\,.
\ee
For $_r\t{\phi}_s$ series the corresponding property, valid independent 
of the factor $\lb (-1)^n (Q/P)^{n(n-1)/2} \rb^{1+s-r}$, is  
\bea
\lim_{a_{rq} \rightarrow \infty}\,_r\t{\phi}_s  
\lb z/a_{rq} \rb & = & {}_r\t{\phi}_s \lb \t{a}_1, \dots, 
\t{a}_{r-1},(0,1);\t{b}_1, \dots, \t{b}_s; \t{q},z \rb\,, \\ 
  &   &  \nn \\ 
\lim_{a_{rp} \rightarrow \infty}\,_r\t{\phi}_s  
\lb z/a_{rp} \rb & = & {}_r\t{\phi}_s \lb \t{a}_1, \dots, 
\t{a}_{r-1},(1,0);\t{b}_1, \dots, \t{b}_s; \t{q},z \rb\,. 
\eea
At times, the parameters $\t{a}$, $\t{b}$, etc., and $\t{q}$ may be 
explicitly indicated in the formulae and equations as $\lb a_p,a_q 
\rb$, $\lb b_p,b_q \rb$, etc., and $\lb P,Q \rb$, respectively, and 
such notations should be clear from the context.  \\

With the notation
\be
\t{q} ^x = \lb P^x,Q^x \rb \,,
\ee
choosing  
\be
\t{a}_1 = \t{q} ^{\a _1},\ \ \dots \,,\ \ 
\t{a}_r = \t{q} ^{\a _r}\,, \quad 
\t{b}_1 = \t{q} ^{\b _1},\ \ \dots \,,\ \ 
\t{b}_s = \t{q} ^{\b _s}\,,   
\ee
and defining
\bea
[x]_{P,Q} & = & \frac{P^x-Q^x}{P-Q}\,, \label{pqn} \\
   &   & \nn \\ 
\lb [x]_{P,Q} \rb _n & = & \lc 
\ba{ll}
1, & n = 0, \cr
[x]_{P,Q} [x+1]_{P,Q} [x+2]_{P,Q} \dots & \cr 
\qquad \dots\,[x+n-1]_{P,Q}\,, & n = 1,2,\dots ,
\ea \right. \\ 
   &   & \nn \\ 
{[n]_{P,Q}}! & = & \lc 
\ba{ll}
1, & n = 0, \cr 
[n]_{P,Q} [n-1]_{P,Q} [n-2]_{P,Q} \dots & \cr 
\qquad \dots\,[2]_{P,Q} [1]_{P,Q}\,, & n = 1,2,\dots\,,
\ea \right.
\eea
it is seen that $_r\tilde{\phi}_s$ becomes 
\bea
   &   &  _r\tilde{\phi}_s \lb \t{q} ^{\a _1},\t{q} ^{\a_ 2},\dots , 
\t{q} ^{\a_ r};\t{q} ^{\b _1},\t{q} ^{\b _2}, \dots , \t{q} ^{\b _s};
\t{q},z \rb \nn \\ 
   &   & \ = \s \,\frac{ \lb \lsb \a _1 \rsb _{P,Q} \rb _n 
\lb \lsb \a _2 \rsb _{P,Q} \rb _n \dots \lb \lsb \a _r \rsb _{P,Q} \rb 
_n }
{ \lb \lsb \b _1 \rsb _{P,Q} \rb _n \lb \lsb \b _2 \rsb _{P,Q} \rb _n 
\dots \lb \lsb \b _s \rsb _{P,Q} \rb _n }
\frac{ \lb (-1)^n (Q/P)^{n(n-1)/2} \rb ^{(1+s-r)} } 
{ (P-Q)^{n(1+s-r)} [n]_{P,Q}! }\,z^n\,. \nn \\ 
   &   & 
\eea 

Note that 
\be
[x]_{P,Q} = P^{x-1} \lb \frac{ 1-\rho^x }{ 1-\rho } \rb =
 P^{x-1} [x]_\rho \,, \quad {\rm with}\ \ \rho = Q/P\,.
\ee 
Thus we can define $_r\tilde{\phi}_s$ also as 
\bea
   &   & _r\tilde{\phi}_s \lb \t{q} ^{\a _1},\t{q} ^{\a_ 2},\dots , 
\t{q} ^{\a_ r};\t{q} ^{\b _1},\t{q} ^{\b _2}, \dots , 
\t{q} ^{\b _s};\t{q},z \rb \nn \\
   &   & \ \ = \s \,\frac{ \lb \lsb \a _1 \rsb _\rho \rb _n 
\lb \lsb \a _2 \rsb _\rho \rb _n \dots \lb \lsb \a _r \rsb _\rho \rb 
_n }
{ \lb \lsb \b _1 \rsb _\rho \rb _n \lb \lsb \b _2 \rsb _\rho \rb _n 
\dots \lb \lsb \b _s \rsb _\rho \rb _n } 
\frac{ \lb (-1)^n (\rho/P)^{n(n-1)/2} \rb ^{1+s-r} }
{ (1-\rho )^{n(1+s-r)} [n]_\rho ! }\,\omega^n\,, \nn \\
   &   & \qquad \qquad \qquad {\rm with}\ \ 
\omega = P^{\lb \sum_{i=1}^r\,\a _i - \sum_{i=1}^s\,\b _i -1 \rb } z\,.  
\eea 
When $(P,Q) = (1,q)$ it is seen that $_r\tilde{\phi}_s \lb 
\t{q} ^{\a _1},\t{q} ^{\a_ 2},\dots , \t{q} ^{\a_ r};\t{q} ^{\b 
_1},\t{q} ^{\b _2}, \dots , \t{q} ^{\b _s};\t{q},z \rb$ reduces to  
$_r\phi_s \lb q^{\a _1},q^{\a_ 2},\dots , q^{\a_ r};q^{\b _1},
q^{\b _2}, \dots , q^{\b _s};q,z \rb$. \\

When $s = r-1$, we have  
\bea
   &   & _r\tilde{\phi}_{r-1} \lb \t{q} ^{\a _1},\t{q} ^{\a_2},\dots , 
\t{q} ^{\a_ r};\t{q} ^{\b _1},\t{q} ^{\b _2}, \dots , 
\t{q} ^{\b _s};\t{q},z \rb \nn \\
   &   & \ \ = {}_r\phi_{r-1} \lb \rho ^{\a _1},\rho ^{\a_2},\dots , 
\rho ^{\a_ r};\rho ^{\b _1},\rho ^{\b _2}, \dots , 
\rho ^{\b _s};\rho , \omega \rb \,. 
\eea
This $_r\tilde{\phi}_{r-1}$ is exactly the $(p,q)$-hypergeometric 
series studied in detail by Burban and Klimyk~\cite{BK}, except for 
the difference in notations and choice of parameters; our  
$\tilde{\phi}$, $r$, $P$, $Q$ and $\rho$ correspond, respectively, 
to their $\Psi$, $A$, $q^{-1/2}$, $p^{1/2}$ and $r$. 

\bigskip

\noindent{\bf 4. Genesis of $(P,Q)$-Analysis} 
\renewcommand{\theequation}{4.{\arabic{equation}}}
\setcounter{equation}{0}

\smallskip

\noindent Let me now recall briefly the genesis of the $(P,Q)$-basic 
number and $(P,Q)$-analysis.  In 1991, Ranabir Chakrabarti and I   
introduced in~\cite{CJ} the $(p,q)$-oscillator algebra 
\be
aa^\dagger - qa^\dagger a = p^{-N}\,, \quad 
[N,a] = -a\,, \quad \lsb N,a^\dagger \rsb = a^\dagger\,, 
\label{pqo}
\ee
generalizing/unifying several forms of $q$-oscillators well 
known in the earlier physics literature related to quantum groups.  
We related the algebra~(\ref{pqo}) to the realization of a 
$(p,q)$-deformed angular momentum algebra
\be
\lsb J_0 , J_\pm \rsb = \pm J_\pm\,, \quad 
J_+J_- - pq^{-1} J_-J_+ = \frac{p^{-2J_0}-q^{2J_0}}{p^{-1}-q}\,.
\label{pqam}
\ee  
Note that the algebra~(\ref{pqo}) is satisfied when 
\be
a^\dagger a = \frac{p^{-N}-q^N}{p^{-1}-q}\,, \quad 
aa^\dagger = \frac{p^{-(N+1)}-q^{N+1}}{p^{-1}-q}\,. 
\ee 
Thus we were led to study the $\lb p^{-1},q \rb$-basic number 
$\lb p^{-n}-q^n \rb / \lb p^{-1}-q \rb$, or $[n]_{p^{-1},q}$ in 
the notation of~(\ref{pqn}), which is a solution for $F_n$ obeying  
the Fibonacci relation 
$F_{n+1} - \lb p^{-1}+q \rb F_n-p^{-1}qF_{n-1} = 0$ for $n \geq 1$, 
with $F_1 = 1$ and $F_0 = 0$.  Further, study of the Bargmann-Fock 
realization and coherent states of the $(p,q)$-oscillator led us to 
define the corresponding $\lb p^{-1},q \rb$-deformation of 
differentiation, integration (for monomials), and exponential 
($_0\t{\phi}_0 \lb -;-;\lb p^{-1},q \rb , \lb p^{-1}-q \rb z \rb$ 
in the notation of the present work, if the factor 
$\lb (-1)^n \lb (Q/P) \rb^{n(n-1)/2} \rb^{1+s-r}$ is dropped from the 
definition~(\ref{tphi}) for $_r\t{\phi}_s$).  With the addition of a 
central element, the $(p,q)$-angular momentum algebra~(\ref{pqam}) 
can be turned into a genuine two-parameter Hopf algebra 
$U_{p,q}(gl(2))$.  Representation theory of $U_{p,q}(gl(2))$ is 
very similar to that of $U_q(sl(2))$ and can be used to construct 
all the finite-dimensional representations of the Hopf algebra 
$Fun_{p,q}(GL(2))$ dual to $U_{p,q}(gl(2))$ (see~\cite{JVdJ}). \\

Around the same time in 1991, independently, there appeared two other 
similar, but very much less detailed, works in the context of quantum 
groups~: Brodimas, Jannussis and Mignani also introduced the 
$(p,q)$-oscillator algebra~(\ref{pqo}) and defined the two-parameter 
deformed derivative~\cite{BJM}.  Arik {\em et al.} also introduced the 
two-parameter oscillator algebra~(\ref{pqo}) calling it the Fibonacci 
oscillator~\cite{Aet}. \\

It is really a surprising coincidence that in the same year 1991, 
without any connection to the quantum group related mathematics/physics 
literature, there appeared a paper~\cite{WW} in which the 
$(p,q)$-basic number, defined by $\lb p^n-q^n \rb / \lb p-q \rb$, was 
introduced while generalizing the Sterling numbers, motivated by 
certain combinatorial problems.  This work of Wachs and White~\cite{WW} 
was brought to the notice of physicists by Katriel and Kibler~\cite{KK} 
who also defined the $(p,q)$-binomial coefficients and derived a
$(p,q)$-binomial theorem while discussing normal ordering for deformed
boson operators obeying algebra of the type~(\ref{pqo}).  Smirnov and 
Wehrhahn~\cite{SW} gave an operator version of such a $(p,q)$-binomial 
theorem giving an expression for the expansion of 
$\lb q^{J_0(1)}J_\pm(2) + J_\pm(1)p^{-J_0(2)} \rb^l$ in terms of the 
$(p,q)$-binomial coefficients, where $\left\{ J_0(1), J_\pm(1) \right\}$ 
and $\left\{ J_0(2), J_\pm(2) \right\}$ are the generators of two 
commuting $(p,q)$-angular momentum algebras~(\ref{pqam}). \\

Before closing, I like to acknowledge that the idea of introducing 
$\lb \t{x};\t{q} \rb _n$, as defined in ~(\ref{xqn}), is derived from the 
definitions
\be
(\lambda ; x)^{(l)} = (\lambda +x)(p\lambda + qx) \lb p^2\lambda 
+q^2x \rb \dots \lb p^{l-1}\lambda + q^{l-1}x \rb\,,
\ee 
and 
\be
\lsb p^\mu,p^\nu;p,q \rsb _n = \lb \frac{1}{p^\mu} - q^\nu \rb 
\lb \frac{1}{p^{\mu+1}} - q^{\nu+1} \rb \dots 
\lb \frac{1}{p^{\mu+n-1}} - q^{\nu+n-1} \rb\,, 
\label{piqn}
\ee 
occurring in~\cite{KK} and~\cite{FLV}, respectively, in related, but 
different, contexts.  In our notation, $(\lambda ; x)^{(l)}$ 
of~\cite{KK} is $((\lambda, -x);(p,q))_l$ and $\lsb p^\mu,p^\nu;p,q 
\rsb _n$ of~\cite{FLV} is $\lb \lb p^{-\mu},q^\nu \rb ; \lb 
p^{-1},q \rb \rb_n$.  In~\cite{KK} the $(p,q)$-binomial coefficients 
have been defined and a $(p,q)$-analogue of the binomial theorem 
$(a+b)^n = \sum_{k=0}^n\,\left( \begin{array}{c}
n \\
k 
\end{array} \right)a^kb^{n-k}\,,\ n = 0,1,2,\ldots\,,$ has been 
obtained.  In~\cite{FLV} the $(p,q)$-oscillator algebra~(\ref{pqo}) 
and $\lsb p^\mu,p^\nu;p,q \rsb _n$ are related to bibasic hypergeometric 
functions~\cite{AV} regarding $p$ and $q$ as two different bases.  
Further in this connection, let me note that the construction of 
$(P,Q)$-analogue of multibasic hypergeometric series is straightforward~: 
one has to choose multiple $\t{q}$-doublets ($\t{q}_1 = \lb P_1,Q_1 
\rb ,\t{q}_2 = \lb P_2,Q_2 \rb , \dots$) and assign them partially 
to the numerator parameters ($\t{a}$-doublets) and the denominator 
parameters ($\t{b}$-doublets), thus simply extending the procedure 
adopted in constructing the usual multibasic hypergeometric series 
(see~\cite{AV} and, {\em e.g.},~\cite{GR}). 

\bigskip 

\noindent{\bf 5. Conclusion} 
\renewcommand{\theequation}{5.{\arabic{equation}}}
\setcounter{equation}{0}

\smallskip

\noindent To conclude, let me just give one example which, I hope, will 
convince the reader that this work indeed promises to lead to an 
interesting subject: $(P,Q)$-special functions.  It is straightforward 
to obtain a general $(P,Q)$-binomial theorem~:    
\bea
   &   & _1\t{\phi}_0 \lb \lb a_p,a_q \rb ;-; (P,Q),z \rb = 
\s \,\frac{ \lb \lb a_p,a_q \rb ; \lb P,Q 
\rb \rb _n } { \lb \lb P,Q \rb ; \lb P,Q \rb \rb _n }\,z^n = 
\frac{ \lb \lb P,a_q z \rb ; \lb P,Q \rb \rb _\infty } 
{ \lb \lb P,a_p z \rb ; \lb P,Q \rb \rb _\infty }\,. \nn \\ 
   &   & 
\label{binomi}
\eea 
In the special case when $(P,Q) = \lb q^{-1/2},p^{1/2} \rb$ and 
$\lb a_p,a_q \rb = \lb q^{-a/2},p^{a/2} \rb$ equation~(\ref{binomi})
reduces to the $(p,q)$-binomial theorem of Burban and Klimyk~\cite{BK}. 
It is easily seen that this $(P,Q)$-binomial theorem~(\ref{binomi}) has 
interesting consequences.  The product 
$\prod_{i=1}^n\,_1\t{\phi}_0 \lb \lb a_{ip},a_{iq} \rb ;-; (P,Q),z 
\rb$ is invariant under the group of independent permutations of 
the $p$-components $\lb a_{1p},a_{2p},\dots ,a_{np} \rb$ and the 
$q$-components $\lb a_{1q},a_{2q},\dots ,a_{nq} \rb$.  Note 
that this product has value $1$ if the $n$-tuple of $p$-components 
$\lb a_{1p},a_{2p},\dots ,a_{np} \rb$ is related to the $n$-tuple of 
$q$-components $\lb a_{1q},a_{2q},\dots ,a_{nq} \rb$ by mere permutation. 
A special case is the relation $_1\t{\phi}_0 \lb (1,0);-; (1,q),z \rb 
{}_1\t{\phi}_0 \lb (0,1);-; (1,q),z \rb = 1$ which is the well known 
identity $e_q(z)E_q(-z) = 1$ for the two canonical $q$-exponentials. 

\bigskip

\noindent {\bf Acknowledgements}~: This article is based mainly on  
the talk I gave at the Workshop on Special Functions \& Differential 
Equations, Chennai, 1997 (WSSF97); it also contains a few new 
results which were not part of my talk but sprang up naturally 
while writing up this contribution to the Proceedings of the 
Workshop and I could not resist the temptation to include.  I am 
very grateful to Prof. K. Srinivasa Rao, particularly with regard to 
the present work, besides for many other things~: whatever little 
fascination I have got for the general theory of hypergeometric 
series has been induced by him.  Actually, I wish to confess that 
my interest in obtaining a $(p,q)$-generalization of the 
$q$-hypergeometric series has its origin in the period around July 
1991 when I was attracted to certain ideas of Prof. Srinivasa Rao 
in this direction which have at last had a definitive effect on me 
recently.  I am thankful to Prof. J. Van der Jeugt for a copy of the 
paper of Burban and Klimyk.  This research was partly supported by 
the EEC (contract No. CI1$^*$-CT92-0101).

\end{document}